\documentclass[12pt]{article}
\usepackage{amssymb}
\def\vs{\vspace{5mm}}
\def\D{{\mathbb D}}
\def\N{{\mathbb N}}
\def\R{{\mathbb R}}
\def\Z{{\mathbb Z}}
\def\sgn{\hbox{sgn}}
\begin{document}
\centerline{\bf Linear natural liftings of forms}

\centerline{\bf to Weil bundles with Weil algebras $\D^r_k$}\vs

\centerline{\bf Jacek D\c ebecki}

\centerline{Instytut Matematyki Uniwersytetu Jagiello\'nskiego}

\centerline{ul. Reymonta 4, 30-059 Krak\'ow, Poland}

\centerline{e-mail: debecki@im.uj.edu.pl}\vs

{\bf Abstract.} We give an explicit description and calculate the
dimension of the vector space of linear natural liftings of
$p$-forms on $n$-dimensional manifolds $M$ to $q$-forms on
$T^{\D^r_k}M$, where $\D^r_k$ is the Weil algebra of $r$-jets at
$0$ of smooth functions $\R^k\longrightarrow\R$, for all
non-negative integers $n$, $p$, $q$, $r$, $k$ except the case
$p=n$ and $q=0$.\vs



Let $A$ be a Weil algebra and $T^A$ the Weil functor corresponding
to $A$ (see [\ref{KM}] or [\ref{KMS}]). Let us denote by
${\mit\Omega}^pM$ the vector space of $p$-forms on a manifold $M.$
A linear natural lifting of $p$-forms to $q$-forms on $T^A$ is a
family of linear maps
$L_M:{\mit\Omega}^pM\longrightarrow{\mit\Omega}^q(T^AM)$ indexed
by $n$-dimensional manifolds and satisfying for all such manifolds
$M$, $N$, every embedding $f:M\longrightarrow N$ and every
$\omega\in{\mit\Omega}^pN$ the condition
$L_M(f^*\omega)=(T^Af)^*(L_N(\omega))$.

In [\ref D] we have given a classification of liftings of this
kind for all non-negative integers $n$, $p$ and $q$ except the
case $p=n$ and $q=0$. There we have established an isomorphism
between the vector space of such liftings and the vector space in
the table below for proper $n$, $p$ and $q$.\vs

\centerline{\begin{tabular}{|c|c|c|c|}\hline
&$0\le p\le n-1$&$p=n$&$n+1\le p$\\\hline
$q=0$&$A_{p-q}$&&$\{0\}$\\\hline
$1\le q\le p$&$A_{p-q}\oplus A_{p-q+1}$&$A_{p-q}$&$\{0\}$\\\hline
$q=p+1$&$A_{p-q+1}$&$\{0\}$&$\{0\}$\\\hline
$p+2\le q$&$\{0\}$&$\{0\}$&$\{0\}$\\\hline
\end{tabular}}\vs

\noindent Here $A_s$ for the Weil algebra $A$ inducing $T^A$ and a non-negative
integer $s$ is the vector space of skew-symmetric $s$-linear maps
$F:A\times\ldots\times A\longrightarrow A^*$, where $A^*$ denotes the
vector space of linear functions $A\longrightarrow\R$, satisfying
\begin{eqnarray}\label1
\lefteqn{F(a_1,\ldots,a_{t-1},bc,a_{t+1},\ldots,a_s)(d)=}\\
\nonumber&&F(a_1,\ldots,a_{t-1},b,a_{t+1},\ldots,a_s)(cd)+F(a_1,\ldots,a_
{t-1},c,a_{t+1},\ldots,a_s)(bd)
\end{eqnarray}
for every $t\in\{1,\ldots,s\}$ and all
$a_1,\ldots,a_{t-1},a_{t+1},\ldots,a_s,b,c,d\in A$.

Unfortunately, the vector spaces $A_s$ are a bit abstract and we cannot
find out the dimension of $A_s$ for every Weil algebra $A$ and every
non-negative integer $s$. This work is devoted to the study of a special
case, namely if $A$ is the algebra $\D^r_k$ of $r$-jets at $0$ of smooth
functions $\R^k\longrightarrow\R$. We will give an explicit description
of $(\D^r_k)_s$ and calculate its dimension for all non-negative integers
$r$, $k$, $s$. The importance of the special case we treat is that each
Weil algebra $A$ is a factor algebra of $\D^r_k$ for some $r$, $k$ (see
[\ref K]), so $A_s$ is a subspace of $(\D^r_k)_s$ for each $s$.

Fix non-negative integers $r$, $k$, $s$. We will denote by $x^i$ for
$i\in\{1,\ldots,k\}$ the $r$-jet at $0$ of the function
$\R^k\ni u\longrightarrow u^i\in\R$ and we will write
$x^\alpha=(x^1)^{\alpha^1}\ldots(x^k)^{\alpha^k}$ and
$|\alpha|=\alpha^1+\ldots+\alpha^k$ for each
$\alpha\in\N^k$, where $\N$ stands for the set of non-negative
integers. It is obvious that $x^\varepsilon$ for $\varepsilon\in\N^k$
such that $|\varepsilon|\le r$ form a basis of the vector space $\D^r_k$
and
$$x^\zeta x^\eta=\cases{x^{\zeta+\eta}&if $|\zeta+\eta|\le r$,\cr0
&otherwise}$$
for all $\zeta,\eta \in\N^k$ such that $|\zeta|\le r$, $|\eta|\le r$.

Of course, $A_0=A^*$ for every Weil algebra $A$. Therefore we will be
concerned only with the case $s>0$. If $r=0$ or $k=0$ then $\D^r_k=\R$,
so $(\D^r_k)_0=\R^*$ and it is a simple matter to see that if $s>0$ then
$(\D^r_k)_s=\{0\}$. Therefore we will be concerned only with the case
$r>0$ and $k>0$.

We can now formulate our main result.\vs

{\bf Definition.} Let $Z$ denote the set of
$(i_1,\ldots,i_s,\alpha)\in\{1,\ldots,k\}^s\times\N^k$ with the
properties that $i_1<\ldots<i_s$ and either $|\alpha|<r$ or
$|\alpha|=r$ and $i_s<\max\{l\in\{1,\ldots,k\}:\alpha^l>0\}$.\vs

{\bf Theorem.} {\it The map $I:(\D^r_k)_s\longrightarrow\R^Z$ given by
$$I(F)(i_1,\ldots,i_s,\alpha)=F(x^{i_1},\ldots,x^{i_s})(x^\alpha)$$
for every $F\in(\D^r_k)_s$ and every $(i_1,\ldots,i_s,\alpha)\in
Z$ is an isomorphism of vector spaces.\vs

Proof.} The theorem will be proved by showing that for each $C\in\R^Z$
there is a unique $F\in(\D^r_k)_s$ such that
\begin{eqnarray}\label2
F(x^{i_1},\ldots,x^{i_s})(x^\alpha)=C(i_1,\ldots,i_s,\alpha)
\end{eqnarray}
for every $(i_1,\ldots,i_s,\alpha)\in Z$. Fix $C\in\R^Z$. Our
construction of $F$ will be divided into six steps.

{\it Step 1.} We define $F(x^{i_1},\ldots,x^{i_s})(x^\alpha)$ for
$(i_1,\ldots,i_s,\alpha)\in Z$ by (\ref2).

{\it Step 2.} We define $F(x^{i_1},\ldots,x^{i_s})(x^\alpha)$ for
$(i_1,\ldots,i_s,\alpha)\in\{1,\ldots,k\}^s\times\N^k$ such that either
$|\alpha|<r$ or $|\alpha|=r$ and
$\max\{i_1,\ldots,i_s\}<\max\{l\in\{1,\ldots,k\}:\alpha^l>0\}$.

Since $F$ should be skew-symmetric, we put either
$$F(x^{i_1},\ldots,x^{i_s})(x^\alpha)=\sgn\sigma F(x^{i_{\sigma(1)}},
\ldots,x^{i_{\sigma(s)}})(x^\alpha)$$
if there is a permutation $\sigma$ of $\{1,\ldots,s\}$ such that
$i_{\sigma(1)}<\ldots< i_{\sigma(s)}$ (note that there is at most one
$\sigma$ with this property) or
$F(x^{i_1},\ldots,x^{i_s})(x^\alpha)=0$ otherwise.

{\it Step 3.} We define $F(x^{i_1},\ldots,x^{i_s})(x^\alpha)$ for
$(i_1,\ldots,i_s,\alpha)\in\{1,\ldots,k\}^s\times\N^k$ such that
$i_1<\ldots<i_s$ and $|\alpha|\le r$, but
$(i_1,\ldots,i_s,\alpha)\notin Z$.

If $G\in(\D^r_k)_s$, $t\in\{1,\ldots,s\}$ and
$\gamma_1,\ldots,\gamma_s,\delta\in\N^k$ then, by induction on
$|\gamma_t|$, (\ref1) leads easily to
\begin{eqnarray}\label3
\lefteqn{G(x^{\gamma_1},\ldots,x^{\gamma_s})(x^\delta)=}\\
\nonumber&&\sum_{j\in\{l\in\{1,\ldots,k\}:\gamma^l_t>0\}}\gamma^j_tG(x^{
\gamma_1},\ldots,x^{\gamma_{t-1}},x^j,x^{\gamma_{t+1}},\ldots,x^{\gamma_s
})(x^{\gamma_t-e_j+\delta}),
\end{eqnarray}
where $e_1,\ldots,e_k$ stand for the standard basis of the module $\Z^k$.

The condition $(i_1,\ldots,i_s,\alpha)\notin Z$ means that $|\alpha|=r$
and $i_s\ge\max\{l\in\{1,\ldots,k\}:\alpha^l>0\}$. Taking $t=s$,
$\gamma_1=e_{i_1},\ldots,\gamma_{s-1}=e_{i_{s-1}}$,
$\gamma_s=\alpha+e_{i_s}$ and $\delta=0$ in (\ref3) we see that $F$
should satisfy
\begin{eqnarray}\label4
\lefteqn{0=(\alpha^{i_s}+1)F(x^{i_1},\ldots,x^{i_s})(x^\alpha)+}\\
\nonumber&&\hspace{103pt}\sum_{j\in\{l\in\{1,\ldots,k\}:\alpha^l>0\}
\setminus\{i_s\}}\alpha^jF(x^{i_1},\ldots,x^{i_{s-1}},x^j)(x^{\alpha+e_{i
_s}-e_j}),
\end{eqnarray}
since $x^{\alpha+e_{i_s}}=0$ and $F$ should be $s$-linear. But
$F(x^{i_1},\ldots,x^{i_{s-1}},x^j)(x^{\alpha+e_{i_s}-e_j})$ for every
$j\in\{l\in\{1,\ldots,k\}:\alpha^l>0\}\setminus\{i_s\}$
has already been defined, as
$\max\{i_1,\ldots,i_{s-1},j\}<i_s=\max\{l\in\{1,\ldots,k\}:(\alpha+e_{i_s
}-e_j)^l>0\}$. Therefore we put
\begin{eqnarray}\label5
\lefteqn{F(x^{i_1},\ldots,x^{i_s})(x^\alpha)=}\\
\nonumber&&-{1\over\alpha^{i_s}+1}\sum_{j\in\{l\in\{1,\ldots,k\}:\alpha^l
>0\}\setminus\{i_s\}}\alpha^jF(x^{i_1},\ldots,x^{i_{s-1}},x^j)(x^{\alpha+
e_{i_s}-e_j}).
\end{eqnarray}

{\it Step 4.} We define $F(x^{i_1},\ldots,x^{i_s})(x^\alpha)$ for
$(i_1,\ldots,i_s,\alpha)\in\{1,\ldots,k\}^s\times\N^k$ such that
$|\alpha|\le r$.

This goes in the same way as step 2.

{\it Step 5.} We define $F(x^{\gamma_1},\ldots,x^{\gamma_s})(x^\delta)$
for $\gamma_1,\ldots,\gamma_s,\delta\in\N^k$ such that
$|\gamma_1|\le r,\ldots,|\gamma_s|\le r,|\delta|\le r$.

Since $F$ should satisfy (\ref3) and take linear values, we put either
\begin{eqnarray}\label6
\lefteqn{F(x^{\gamma_1},\ldots,x^{\gamma_s})(x^\delta)=\sum_{j_1\in\{l\in
\{1,\ldots,k\}:\gamma^l_1>0\}}\ldots\sum_{j_s\in\{l\in\{1,\ldots,k\}:
\gamma^l_s>0\}}}\\
\nonumber&&\hspace{148pt}\gamma^{j_1}_1\ldots\gamma^{j_s}_sF(x^{j_1},
\ldots,x^{j_s})(x^{\gamma_1-e_{j_1}+\ldots+\gamma_s-e_{j_s}+\delta})
\end{eqnarray}
if $|\gamma_1+\ldots+\gamma_s+\delta|\le r+s$ or
$F(x^{\gamma_1},\ldots,x^{\gamma_s})(x^\delta)=0$ otherwise.

{\it Step 6.} We complete our construction easily, because
$x^\varepsilon$ for $\varepsilon\in\N^k$ such that $|\varepsilon|\le r$
form a basis of the vector space $\D^r_k$ and $F$ should be $s$-linear
with linear values.

Thus we have proved the uniqueness of $F$. By step 1, the map $F$ we have
constructed satisfies (\ref2). By step 6, it is $s$-linear with linear
values. By steps 4, 5, 6, it is easily seen to be skew-symmetric. What is
left is to prove that it satisfies (\ref1).

We begin by showing the crucial fact that for all
$g_1,\ldots,g_{s-1}\in\{1,\ldots,k\}$ such that $g_1<\ldots<g_{s-1}$ and
every $\varepsilon\in\N^k$ such that $|\varepsilon|=r+1$
\begin{equation}\label7
\sum_{h\in\{l\in\{1,\ldots,k\}:\varepsilon^l>0\}}\varepsilon^hF(x^{g_1},
\ldots,x^{g_{s-1}},x^h)(x^{\varepsilon-e_h})=0.
\end{equation}

If either $s=1$ or $s>1$ and $g_{s-1}<m$, where
$m=\max\{l\in\{1,\ldots,k\}:\varepsilon^l>0\}$ which implies
$m\ge\max\{l\in\{1,\ldots,k\}:(\varepsilon-e_m)^l>0\}$, then (\ref7) is
nothing but (\ref4) with $i_1=g_1,\ldots,i_{s-1}=g_{s-1}$, $i_s=m$ and
$\alpha=\varepsilon-e_m$. So (\ref7) holds, as (\ref4) is equivalent to
(\ref5) which holds.

If $s>1$ and $g_{s-1}\ge m$, then
$\max\{g_1,\ldots,g_{s-1},h\}=g_{s-1}\ge\max\{l\in\{1,\ldots,k\}:(
\varepsilon-e_h)^l>0\}$ for every
$h\in\{l\in\{1,\ldots,k\}:\varepsilon^l>0\}$. Hence if
$h\notin\{g_1,\ldots,g_{s-1}\}$ then the skew-symmetry of $F$ and (\ref5)
with $\{i_1,\ldots,i_s\}=\{g_1,\ldots,g_{s-1},h\}$ (which implies
$i_s=g_{s-1}$) and $\alpha=\varepsilon-e_h$ give
\begin{eqnarray}\label8
\lefteqn{F(x^{g_1},\ldots,x^{g_{s-1}},x^h)(x^{\varepsilon-e_h})=-{1\over
\varepsilon^{g_{s-1}}+1}\sum_{j\in\{l\in\{1,\ldots,k\}:\varepsilon^l>0\}
\setminus\{g_{s-1},h\}}}\\
\nonumber&&\hspace{160pt}\varepsilon^jF(x^{g_1},\ldots,x^{g_{s-2}},x^j,x^
h)(x^{\varepsilon-e_h+e_{g_{s-1}}-e_j}),
\end{eqnarray}
because if $\varepsilon^h>1$ then
$F(x^{g_1},\ldots,x^{g_{s-2}},x^h,x^h)(x^{\varepsilon-e_h+e_{g_{s-1}}-e_h
})=0$, by the skew-symmetry of $F$. Substituting (\ref8) into (\ref7) and
omitting the terms which vanish on account of the skew-symmetry of $F$ we
see that the left hand side of (\ref7) equals
\begin{eqnarray*}
\lefteqn{-{1\over\varepsilon^{g_{s-1}}+1}\sum_{h\in\{l\in\{1,\ldots,k\}:
\varepsilon^l>0\}\setminus\{g_1,\ldots,g_{s-1}\}}\sum_{j\in\{l\in\{1,
\ldots,k\}:\varepsilon^l>0\}\setminus\{g_1,\ldots,g_{s-1},h\}}}\\
&&\hspace{149pt}\varepsilon^h\varepsilon^jF(x^{g_1},\ldots,x^{g_{s-2}},x^
j,x^h)(x^{\varepsilon-e_h+e_{g_{s-1}}-e_j}).
\end{eqnarray*}
But, using the skew-symmetry of $F$ again, we have
\begin{eqnarray*}
\lefteqn{\varepsilon^j\varepsilon^hF(x^{g_1},\ldots,x^{g_{s-2}},x^h,x^j)(
x^{\varepsilon-e_j+e_{g_{s-1}}-e_h})=}\\
&&-\varepsilon^h\varepsilon^jF(x^{g_1},\ldots,x^{g_{s-2}},x^j,x^h)(x^{
\varepsilon-e_h+e_{g_{s-1}}-e_j}).
\end{eqnarray*}
Therefore the left hand side of (\ref7) equals $0$. This establishes
(\ref7).

We recall that our aim is to show (\ref1) for $F$ we have
constructed. Since $F$ is skew-symmetric, it suffices to prove
(\ref1) only for $t=s$. This will be proved as soon as we can show
that
\begin{eqnarray}\label9
\lefteqn{F(x^{\alpha_1},\ldots,x^{\alpha_{s-1}},x^{\beta+\gamma})(x^
\delta)=}\\
\nonumber&&F(x^{\alpha_1},\ldots,x^{\alpha_{s-1}},x^\beta)(x^{\gamma+
\delta})+F(x^{\alpha_1},\ldots,x^{\alpha_{s-1}},x^\gamma)(x^{\beta+\delta
})
\end{eqnarray}
for $\alpha_1,\ldots,\alpha_{s-1},\beta,\gamma,\delta\in\N^k$ such that
$|\alpha_1|\le r,\ldots,|\alpha_{s-1}|\le r$, $|\beta|\le r$,
$|\gamma|\le r$, $|\delta|\le r$, because both the sides of (\ref1)
are $(s+2)$-linear with respect to $a_1,\ldots,a_{s-1}$, $b$, $c$, $d$
and $x^\varepsilon$ for $\varepsilon\in\N^k$ such that
$|\varepsilon|\le r$ form a basis of the vector space $\D^r_k$.

We now observe that (\ref9) holds in four special cases.

{\it Case 1.} $\beta=0$ or $\gamma=0$. Then (\ref9) is evident because,
by steps 5 and 6, we have
$F(x^{\alpha_1},\ldots,x^{\alpha_{s-1}},1)(x^{\gamma+\delta})=0$ and
$F(x^{\alpha_1},\ldots,x^{\alpha_{s-1}},1)(x^{\beta+\delta})=0$.

{\it Case 2.} There is $i\in\{1,\ldots,s-1\}$ such that $\alpha_i=0$.
Then (\ref9) is evident, as both the sides of (\ref9) equal $0$, by
steps 5 and 6.

{\it Case 3.}  $|\alpha_1+\ldots+\alpha_{s-1}+\beta+\gamma+\delta|>r+s$.
Then (\ref9) is also evident, as both the sides of (\ref9) equal $0$, by
steps 5 and 6.

{\it Case 4.}
$|\alpha_1+\ldots+\alpha_{s-1}+\beta+\gamma+\delta|\le r+s$,
$|\beta+\gamma|\le r$, $|\beta+\delta|\le r$, $|\gamma+\delta|\le r$.
Then (\ref9) can be easily checked, because we may apply (\ref6) to the
left hand side of (\ref9) as well as to each of two terms of its right
hand side.

Assume that it is none of the above cases. Since it is not case 2,
$|\alpha_1|\ge1,\ldots,|\alpha_{s-1}|\ge1$. Since it is not case 3,
$|\alpha_1+\ldots+\alpha_{s-1}+\beta+\gamma+\delta|\le r+s$. Combining
these yields $|\beta+\gamma+\delta|\le r+1$. If it were true that
$|\beta+\gamma|\le r$, it would also be true that $|\beta+\delta|>r$ or
$|\gamma+\delta|>r$, as it is not case 4, and so that $\gamma=0$ or
$\beta=0$ respectively, contrary  to the fact that it is not case 1.
Therefore $|\beta+\gamma|=r+1$, and so $\delta=0$,
$|\alpha_1|=1,\ldots,|\alpha_{s-1}|=1$.

Summing up, it remains to prove (\ref9) only if
$|\alpha_1|=1,\ldots,|\alpha_{s-1}|=1$, $|\beta|\le r$, $|\gamma|\le r$,
$|\beta+\gamma|=r+1$ and $\delta=0$. Then there are
$g_1,\ldots,g_{s-1}\in\{1,\ldots,k\}$ such that
$\alpha_1=e_{g_1},\ldots,\alpha_1=e_{g_{s-1}}$. Since $F$ is
skew-symmetric, without loss of generality we can assume that
$g_1<\ldots<g_{s-1}$. Moreover, $x^{\beta+\gamma}=0$ and $F$ is
$s$-linear, hence the left hand side of (\ref9) equals $0$. Using (\ref6)
we can rewrite (\ref9) as
\begin{eqnarray*}
\lefteqn{0=\sum_{h\in\{l\in\{1,\ldots,k\}:\beta^l>0\}}\beta^hF(x^{g_1},
\ldots,x^{g_{s-1}},x^h)(x^{\beta-e_h+\gamma})+}\\
&&\hspace{124pt}\sum_{h\in\{l\in\{1,\ldots,k\}:\gamma^l>0\}}\gamma^hF(x^{
g_1},\ldots,x^{g_{s-1}},x^h)(x^{\gamma-e_h+\beta}),
\end{eqnarray*}
which is nothing but (\ref7) with $\varepsilon=\beta+\gamma$. This
completes the proof of the theorem.\vs

{\bf Corollary.} {\it The dimension of the vector space $(\D^r_k)_s$
equals
$${r+s-1\choose s}{r+k\choose r+s}.$$\vs

Proof.} We will compute the number of elements of $Z$, which is equal to
the dimension of $(\D^r_k)_s$, by the theorem.

For each $v\in\{0,\ldots,r-1\}$ the number of
$(i_1,\ldots,i_s,\alpha)\in\{1,\ldots,k\}^s\times\N^k$ such that
$i_1<\ldots<i_s$ and $|\alpha|=v$ equals
$${k\choose s}{v+k-1\choose k-1}.$$
Furthermore, we have
$${k\choose s}\sum_{v=0}^{r-1}{v+k-1\choose k-1}={k\choose s}{r+k-1
\choose k}={r+s-1\choose s}{r+k-1\choose r+s-1}.$$

For each $m\in\{s+1,\ldots,k\}$ the number of
$(i_1,\ldots,i_s,\alpha)\in\{1,\ldots,k\}^s\times\N^k$ such that
$i_1<\ldots<i_s$, $|\alpha|=r$ and
$i_s<m=\max\{l\in\{1,\ldots,k\}:\alpha^l>0\}$ equals
$${m-1\choose s}{r+m-2\choose m-1}.$$
Furthermore, we have
\begin{eqnarray*}
\lefteqn{\sum_{m=s+1}^k{m-1\choose s}{r+m-2\choose m-1}=}\\
&&{r+s-1\choose  s}\sum_{m=s+1}^k{r+m-2\choose r+s-1}={r+s-1\choose  s}{r
+k-1\choose r+s}.
\end{eqnarray*}

Hence the number of elements of $Z$ equals
$${r+s-1\choose s}{r+k-1\choose r+s-1}+{r+s-1\choose s}{r+k-1\choose r+s}
\hspace{-1.8pt}=\hspace{-1.8pt}{r+s-1\choose s}{r+k\choose r+s}.$$
This completes the proof of the corollary.\vs

Note that the corollary is still true if $r=0$ or $k=0$ or $s=0$,
as is easy to check.\vs

\centerline{\bf References}\vs

\begin{enumerate}
\item\label DJ. D\c ebecki: {\it Linear liftings of $p$-forms to
$q$-forms on Weil bundles,} Mo\-natsh. Math. 148 (2006), 101--117

\item\label{KM}G. Kainz, P. Michor: {\it Natural transformations
in differential geometry,} Czech. Mat. J. 37 (112) 1987, 584--607

\item\label KI. Kol\'a\v r: {\it Jet-like approach to Weil
bundles,} Seminar Lecture Notes, Masaryk University, Brno, 2001

\item\label{KMS}I. Kol\'a\v r, P. W. Michor, J. Slov\'ak: {\it
Natural operations in differential geometry,} Springer, Berlin,
1993
\end{enumerate}
\end{document}